\documentclass[11pt]{amsart}

\usepackage[margin=3.5cm]{geometry}
\usepackage{amsfonts}
\usepackage{amssymb}
\usepackage{amsmath}
\usepackage{color}
\usepackage{fancyhdr}
\usepackage{amsthm}
\usepackage{hyperref}
\usepackage{verbatim}
\usepackage{float}
\usepackage[all]{xy}
\usepackage{graphicx}
\usepackage{caption}
\usepackage{subcaption}
\usepackage{mathtools}
\usepackage{tikz-cd}
\usepackage{enumitem}
\usepackage{titlesec}
\usepackage{MnSymbol}
\usepackage{relsize}
\usepackage{mathtools}
\usepackage{mathrsfs}  
\usepackage{lipsum}
\usepackage[utf8]{inputenc}
\usepackage[backend=biber, style=alphabetic]{biblatex}

\defbibheading{bibliography}[References]{%
    \begin{center}
    \fontsize{12}{15}\textbf{#1}
    \end{center}
}

\titleformat{\section}
  {\normalfont\fontsize{12}{15}\bfseries}{\thesection}{1em}{}

\titleformat{\subsection}
  {\normalfont\fontsize{12}{15}\bfseries}{\thesubsection}{1em}{}

\newtheorem{theorem}{Theorem}[section]

\newtheorem{lemma}[theorem]{Lemma}

\newtheorem{corollary}[theorem]{Corollary}

\theoremstyle{definition}

\newtheorem{example}[theorem]{Example}

\newtheorem{remark}[theorem]{Remark}

\graphicspath{ {images/} }

\newcommand{\mb}{\mathbb}

\newcommand{\Q}{\mb{Q}}

\newcommand{\Z}{\mb{Z}}
\newcommand{\C}{\mb{C}}

\newcommand*{\sheafhom}{\mathcal{H}\kern -.5pt om}

\bibliography{refs}

\numberwithin{equation}{section}

\title[Cohomology Rings of Fiber Products of Lefschetz Fibrations]{Cohomology Rings of Fiber Products of Lefschetz Fibrations}
\author{Kuan-Wen Chen}
\address{Department of Mathematics, Columbia University, New York, NY 10027}
\email{kc3396@columbia.edu}

\begin{document}
\begin{abstract}
We prove a Lefschetz-type theorem for the cohomology rings of fiber products of two Lefschetz fibrations with disjoint critical loci. As an application, we compute the Hodge numbers of several examples of higher-dimensional Schoen's Calabi–Yau manifolds. 
\end{abstract}

\maketitle

\section{Introduction}

A Schoen's Calabi-Yau threefold is constructed as the fiber product of two general rational elliptic surfaces over $\mathbb{P}^1$. Namikawa \cite{Namikawa} and Grassi-Morisson \cite{GrMo} showed that if $E_i\xrightarrow{f_i} \mathbb{P}^1$ are two general rational elliptic surfaces, then the Picard group of the Schoen's Calabi-Yau threefold $E_1\times_{\mathbb{P}^1} E_2$ is given by $$\textup{Pic}(E_1\times_{\mathbb{P}^1} E_2)\cong \frac{\textup{Pic}(E_1)\oplus\textup{Pic}(E_2)}{\langle f_1^*\mathcal{O}_{\mathbb{P}^1}(1)\sim f_2^*\mathcal{O}_{\mathbb{P}^1}(1)\rangle}.$$ This description of the Picard group plays a crucial role in the proof of the Morrison–Kawamata–Totaro cone conjecture for Schoen Calabi–Yau threefolds by Grassi-Morisson \cite{GrMo}.

In this paper, we investigate to what extent a similar phenomenon holds for higher-degree cohomology of algebraic varieties arising from such fiber-product constructions. Our goal is to understand the structure of higher-degree cohomology groups in this setting, with potential applications, for instance, to the study of the Hodge conjecture for these varieties. The main result of the paper is a partial generalization to higher degree cohomology in any dimensions: we prove a Lefschetz-type theorem for the cohomology rings of fiber products of two Lefschetz fibrations with disjoint critical loci, using mixed Hodge theory.

\begin{theorem} \label{fib}
Let $\tilde{X}_i\xrightarrow{f_i} \mb{P}^1$ $(i=1,2)$ be two Lefschetz fibrations, each arising from the blowup of the base locus of a Lefschetz pencil contained in an ample linear system on a smooth projective variety, and let $m_i=\textup{dim}_{\C}\tilde{X}_i$.
Assume the critical loci of $f_1$ and $f_2$ are disjoint in $\mathbb{P}^1$. Then the canonical morphism of graded rings $$H^{\bullet}(\tilde{X}_1,\Q)\otimes_{H^{\bullet}(\mb{P}^1,\Q)}H^{\bullet}(\tilde{X}_1,\Q)\rightarrow H^{\bullet}(\tilde{X}_1\times_{\mathbb{P}^1}\tilde{X}_2,\Q)$$ is an isomorphism in degree $\leq m_1+m_2-2$, and injective in degree $ m_1+m_2-1$.
\end{theorem}

A typical class of algebraic varieties arising as fiber products of Lefschetz fibrations is provided by Schoen's Calabi-Yau manifolds, in dimension three and higher. These manifolds are obtained as fiber products of two Lefschetz fibrations over $\mathbb{P}^1$, each coming from the blowup of the base locus of an anticanonical Lefschetz pencil on a Fano variety, and with disjoint critical loci. In dimension three, such Calabi-Yau varieties were first introduced by Schoen \cite{Schoen} in the 1980s and subsequently studied by Namikawa \cite{Namikawa} and Grassi-Morisson \cite{GrMo} in the 1990s. Schoen's Calabi-Yau threefolds provide the first nontrivial example of the Morrison-Kawamata-Totaro cone conjecture, where the nef cone is not rational polyhedral \cite{GrMo}.

These manifolds have recently attracted renewed attention. Suzuki \cite{Suz} studied arithmetic properties of Schoen's Calabi–Yau manifolds, while Gachet–Lin–Wang \cite{GLW} proved the cone conjecture for higher-dimensional Schoen's Calabi–Yau varieties. Schoen's Calabi–Yau manifolds also play a central role in Sano’s construction of non-Kähler Calabi–Yau manifolds with arbitrarily large second Betti number \cite{Sa}. Moreover, the results of the present paper have been applied by the author \cite{thesis} to show that Sano’s  non-Kähler Calabi–Yau manifolds satisfy the $\partial\overline{\partial}$-lemma and to compute the signature of the cup product on each $(p,q)$-component of the middle cohomology of these manifolds.

Let  $f_i:\tilde{X}_i\rightarrow \mb{P}^1$ $(i=1,2)$ be two Lefschetz fibrations as in Theorem \ref{fib}. Write $\tilde{X}_i=\textup{Bl}_{B_i}X_i$, the blowup of a smooth projective variety $X_i$ along the base locus $B_i$ of a Lefschetz pencil in an ample linear system. Then, as a consequence of Theorem \ref{fib}, we obtain the following formula for the Hodge numbers of the fiber product $\tilde{X}_1\times_{\mathbb{P}^1}\tilde{X}_2$ in degrees below the middle:

\begin{corollary} \label{hformula1}
     For $p+q< m_1+m_2-1$, $$h^{p,q}(\tilde{X}_1\times_{\mb{P}^1}\tilde{X_2})=h^{p,q}(\tilde{X}_1\times\tilde{X_2})-K^{p-1,q-1}-h^{p-2,q-2}(B_1\times B_2),$$ where $$K^{p-1,q-1}=\textup{dim}_{\C}\textup{ker}\left(H^{p-1,q-1}(B_1\times X_2)\oplus H^{p-1,q-1}(X_1\times B_2)\rightarrow H^{p-1,q-1}(B_1\times B_2)\right).$$
\end{corollary}

We note that $K^{p-1,q-1}$ is entirely determined by the Hodge numbers of $X_i$ and $B_i$ via the Lefschetz hyperplane theorem and the K\"unneth formula. Hence, Corollary \ref{hformula1} expresses the Hodge numbers of $\tilde{X}_1\times_{\mb{P}^1}\tilde{X_2}$ below the middle degree purely in terms of the Hodge numbers of the varieties $X_i$ and the base loci $B_i$.

As an application of Corollary \ref{hformula1}, combined with the standard computation of the Hirzebruch $\chi_y$ genus, we are able to determine the full Hodge diamond of certain higher-dimensional Schoen's Calabi-Yau manifolds. See Section \ref{HNS} for detailed examples.

The proof of Theorem \ref{fib} is inspired by an idea of Kloosterman \footnote{The proof of the theorem is inspired by Kloosterman's answer on Mathoverflow: \url{https://mathoverflow.net/questions/68625/hodge-diamond-of-a-calabi-yau-fourfold/70052##70052}}. We compare the cohomology of the fiber product $\tilde{X}_1\times_{\mathbb{P}^1}\tilde{X}_2$ with that of a singular ample hypersurface $W$ birational to it, using a Mayer–Vietoris argument. Although $W$ is singular, it is cohomologically simpler, since the Lefschetz hyperplane theorem applies to singular ample hypersurfaces.

As noted by Kloosterman, this approach breaks down in degrees close to the middle, where the Lefschetz hyperplane theorem yields only injectivity rather than an isomorphism. In our situation, we instead exploit the fibration structures and compute the cohomology in degree equal to the dimension minus one via the Leray spectral sequence. While this strategy is standard, its use here is motivated by the approach of de Cataldo–Rapagnetta–Saccà in their computation of the Hodge numbers of O’Grady 10 \cite{dCRS}.

\textbf{Acknowledgments.}
The author would like to thank Professor Robert Friedman for his numerous suggestions, explanations, clarifications, and the time he generously spent editing this paper. The author is also grateful to Professors Radu Laza and Hsueh-Yung Lin for their helpful comments. Thanks are due to Hung Chiang and Alex Xu for valuable discussions.

\section{Preliminaries}

\subsection{Cohomology Rings of the Blowups}

We begin by recalling the following classical result on the cohomology rings of the blowups:

\begin{theorem}(\cite[Theorem 7.31]{Voi1})\label{Hodgeblowup}

 Let $X$ be a compact complex manifold and $Z\subset X$ a compact complex submanifold of codimension $r$. Let $\pi:\tilde{X}\rightarrow X$ denote the blowup of $X$ along $Z$, and let $E$ be the exceptional divisor. Denote by $\iota:Z\rightarrow X$ and $\iota':E\rightarrow \tilde{X}$ the inclusion maps. Then \begin{enumerate}
     \item For each $k$, there is a pushout diagram of cohomology \begin{center}
    \begin{tikzcd}
H^k_Z(X) \arrow[r] \arrow[d] & H^k(X) \arrow[d,"\pi^*"] \\
H^k_{E}(\tilde{X}) \arrow[r] & H^k(\tilde{X}), 
\end{tikzcd}
\end{center}
where the vertical maps are induced by the pullback morphism.
\item Identifying $H^k_Z(X)\cong H^{k-2r}(Z)$ and $H^{k}_E(\tilde{X})\cong H^{k-2}(E)$ via the Thom isomorphism, the diagram becomes: \begin{center}
    \begin{tikzcd}
H^{k-2r}(Z) \arrow[r,"\iota_{!}"] \arrow[d,"\Phi"] & H^k(X) \arrow[d,"\pi^*"] \\
H^{k-2}(E) \arrow[r,"\iota'_{!}"] & H^k(\tilde{X}), 
\end{tikzcd}
\end{center}
 where the horizontal maps are the Gysin morphisms. The map $\Phi$ is described as follows: for $\alpha \in H^{k-2r}(Z)$, $$\Phi(\alpha)=(\pi|_{E}^*\alpha)\cup e(Q),$$ where $Q$ is the universal quotient bundle of $\pi|_E^*N_{Z/X}$ on $E\cong \mathbb{P}(N_{Z/X})$, defined by the short exact sequence $$0\rightarrow \mathcal{O}_E(-1)\rightarrow \pi|_E^*N_{Z/X}\rightarrow Q\rightarrow 0,$$ and $e(Q)$ denotes the Euler class of $Q$.
\item Let $h=c_1(\mathcal{O}_E(1))\in H^2(E)$. If $X$ is K\"ahler, then the pushout diagram induces the isomorphism of Hodge structures (with Tate twists ignored for simplicity): $$H^k(X)\oplus \bigoplus_{i=0}^{r-2}H^{k-2i-2}(Z)\xrightarrow{\pi^*+\sum_i \iota'_!\circ(\cup h^i)\circ\pi|_E^* } H^k(\tilde{X}).$$
 \end{enumerate}
\end{theorem}
\begin{proof}
We only explain why the map $\Phi$ can be described as in part (2) of Theorem \ref{Hodgeblowup}.
    
By \cite[Chapter I, \S 6]{BottTu}, the local cohomology group $H^k_Z(X)=H^k(N_{Z/X}, N_{Z/X}-Z)$ can be identified with the compact vertical cohomology $H^k_{cv}(N_{Z/X})$, defined using the complex of smooth differential forms on $N_{Z/X}$ with compact support in the vertical direction. This identification follows from the fact that the local cohomology $H^k_Z(X)$ is isomorphic to the cohomology of the Thom space of the normal bundle $N_{Z/X}$, which is the one-point compactification of $N_{Z/X}$. 

Under this identification, the Thom isomorphism $H^{k-2r}(Z)\rightarrow H^k_{cv}(N_{Z/X})$ is given by the composition of the cup product by the Thom class $Th(N_{Z/X})\in H^{2r}_{cv}(N_{Z/X})$ and the pullback map. Since $\pi|_E^*N_{Z/X}\cong \mathcal{O}_{E}(-1)\oplus Q$ as smooth vector bundles, it follows from \cite[Chapter I, Proposition 6.19]{BottTu} that $$Th(\pi|_E^*N_{Z/X})=p_1^*Th(\mathcal{O}_{E}(-1))\cup p_2^*Th(Q),$$ where $p_1$ and $p_2$ are the projection maps from $\pi|_E^*N_{Z/X}$ to the respective summands. 

We consider the commutative diagram:  
\begin{center}
\begin{tikzcd}
H^{k-2r}(Z) \arrow[rr,"\cup Th(N_{Z/X})"] \arrow[d,"\Phi"] && H^k_{cv}(N_{Z/X}) \arrow[d,"g*"] \\
H^{k-2}(E) \arrow[rr,"\cup Th(\mathcal{O}_{E}(-1))"] && H^k_{cv}(\mathcal{O}_{E}(-1)),
\end{tikzcd}
\end{center}
where the horizontal maps are the Thom isomorphisms and the right vertical map is induced by the composition $g: \mathcal{O}_{E}(-1)\hookrightarrow \pi|_{E}^*N_{Z/X}\rightarrow N_{Z/X}$. For all $\alpha\in H^{k-2r}(Z)$, we compute $$g^*(\alpha\cup  Th(N_{Z/X}))=\pi^*\alpha \cup (Th(\pi|_E^*N_{Z/X})|_{\mathcal{O}_{E}(-1)})$$$$=\pi|_E^*\alpha \cup (p_2^*Th(Q))|_{\mathcal{O}_{E}(-1)}\cup Th(\mathcal{O}_{E}(-1)).$$
We now examine the diagram:     
\begin{center}
\begin{tikzcd}
\mathcal{O}_{E}(-1) \arrow[r] \arrow[d] & \pi|_{E}^*N_{Z/X}= \mathcal{O}_{E}(-1)\oplus Q\arrow[d] \\
E \arrow[r] & Q,
\end{tikzcd}
\end{center}
where the vertical maps are projection maps, the top horizontal map is the inclusion map of vector bundles, and the bottom horizontal map is the embedding of $E$ as the zero section of $Q$. 

It follows that $(p_2^*Th(Q))|_{\mathcal{O}_{E}(-1)}$ is the pullback of the restriction  $Th(Q)|_E$. By  \cite[Chapter I, Proposition 6.41]{BottTu}, we have  $Th(Q)|_{E}=e(Q)$, the Euler class of $Q$. Therefore, we conclude: $$\Phi(\alpha)=\pi|_E^*\alpha\cup e(Q)$$ as claimed. This completes the proof.

\end{proof}

We will need the following lemma, which describes the restriction and Gysin morphisms between the blowup and the strict transform.

\begin{lemma}\label{gysin}
Let $Z\xrightarrow{\iota_2} Y\xrightarrow{\iota_1} X$ be inclusions of compact complex manifolds of dimension $m-2,\ m-1, \ m$ respectively. Let $\tilde{X}=\textup{Bl}_{Z}X$ , and note that $Y\cong \tilde{Y}:=\textup{Bl}_{Z}Y$ naturally embeds into $\tilde{X}$ via  $\tilde{\iota}_1:\tilde{Y}\hookrightarrow \tilde{X}$. Then, for all $k$, the restriction map $$H^k(\tilde{X})\xrightarrow{\tilde{\iota}_1^*}H^{k}(\tilde{Y})$$ can be identified with the map $$H^{k}(X)\oplus H^{k-2}(Z)\xrightarrow{(\iota_1^*,\ (\iota_2)_{!})} H^k(Y),$$ where $\iota_1^*$ is the restriction map and $(\iota_2)_{!}$ is the Gysin morphism.

Dually, for all $k$, the Gysin morphism $H^{k}(\tilde{Y})\xrightarrow{(\tilde{\iota}_1)_{!}}H^{k+2}(\tilde{X})$ can be identified with $$H^{k}(Y)\xrightarrow{((\iota_1)_{!},\ -\iota_2^*)}H^{k+2}(X)\oplus H^{k}(Z),$$ where $(\iota_1)_{!}$ is the Gysin morphism and $\iota_2^*$ is the restriction map.
\end{lemma}
\begin{proof}
  Let $E\xrightarrow{\iota}\tilde{X}$ denote the exceptional divisor of the blowup morphism $\tilde{X}\xrightarrow{\pi} X$. By Theorem \ref{Hodgeblowup} (1), the cohomology $H^k(\tilde{X})$ fits into the pushout diagram: 

\begin{center}
    \begin{tikzcd}
H^k_Z(X) \arrow[r] \arrow[d] & H^k(X) \arrow[d,"\pi^*"] \\
H^k_{E}(\tilde{X}) \arrow[r] & H^k(\tilde{X}), 
\end{tikzcd}
\end{center}
where the vertical maps are induced by pullback.

Note that the composition $H^k(X)\xrightarrow{\pi^*} H^k(\tilde{X})\xrightarrow{\tilde{\iota}_1^*} H^k(\tilde{Y})$ is simply the restriction $H^k(X)\xrightarrow{\iota_1^*}H^k(Y)$. Meanwhile, the composition $H^k_E(\tilde{X})\rightarrow H^k(\tilde{X})\xrightarrow{\tilde{\iota}_1^*} H^k(\tilde{Y})$ factors through  
$ H^k_{E\cap \tilde{Y}}(\tilde{Y})\rightarrow H^k(\tilde{Y})$. Under the natural identification $E\cap \tilde{Y}\cong Z$, this map corresponds to the Gysin morphism $H^{k-2}(Z)\rightarrow H^k(Y)$, via the Thom isomorphism. Therefore, we can identify the full restriction map as: $$H^k(\tilde{X})\cong H^k(X)\oplus H^{k-2}(Z)\xrightarrow{(\iota_1^*, (\iota_2)_{!})}H^k(Y),$$ as claimed.

The dual statement about the Gysin map follows from Poincar\'e duality and the first statement. The negative sign in the second component arises from the identification $\mathcal{O}_{\tilde{X}}(E)|_{E}=\mathcal{O}_E(-1)$, which contributes a negative sign in the following commutative diagram:

\begin{center}
    \begin{tikzcd}
H_{2m-k-2}(\tilde{X}) \arrow[rrr,"(\pi|_E)_{*}\circ(\bullet . {[E]})"]  &&& H_{2m-k-4}(Z) \\
H^{k+2}(\tilde{X}) \arrow[u,"PD_{\tilde{X}}"] &&& H^k(Z) \arrow[lll,"\iota_!\circ (\pi|_E)^*"'] \arrow[u,"-PD_{Z}"]. 
\end{tikzcd}
\end{center}
 Here, the bottom horizontal arrow is a component of the identification map in Theorem \ref{Hodgeblowup} (3), while the top horizontal arrow is the homological dual in the corresponding degree. The notation $(\bullet .{[E]})$ denotes the intersection product with the class of the exceptional divisor $E$. The vertical maps are given by Poincar\'e duality.

\end{proof}

\subsection{ Lefschetz Fibrations and Vanishing Cohomology}

For general facts about Lefschetz pencils and Lefschetz fibrations, we refer to \cite[Chapter 2]{Voi2} and \cite[ Section C.2.3 ]{PS}.

 Let $X$ be a projective manifold of dimension $m$ and $\{Y_t\}_{t\in \mathbb{P}^1}$ be a Lefschetz pencil on $X$, where $Y_{t}$ is a divisor in some ample linear system $|\mathcal{L}|$ for all $t$. Fix a general member $Y:=Y_{t}$ and let $\iota':Y\hookrightarrow X$ be the inclusion. Let $B$ be the base locus of the pencil. We consider the blowup diagram: 

\begin{center}
    \begin{tikzcd}
E \arrow[r, "\tilde{\iota}"] \arrow[d, "\pi|_E"] & \tilde{X}:=\textup{Bl}_ZX \arrow[d, "\pi"] \\
Z \arrow[r, "\iota"]                      & X,                                         
\end{tikzcd}
\end{center}
where $E$ is the exceptional divisor of the blowup morphism. Note that $$E\cong \mathbb{P}(N_{Z/X})=\mathbb{P}(\iota^*\mathcal{L}\oplus \iota^*\mathcal{L} ).$$ Let $$f:\tilde{X}\hookrightarrow X\times \mathbb{P}^1\rightarrow  \mathbb{P}^1$$ be the associated Lefschetz fibration, and let $\Delta(f)\subset \mathbb{P}^1$ be its critical locus. 

 We work over cohomology with rational coefficients. By the Lefschetz hyperplane theorem, the restriction map $H^k(X)\rightarrow H^{k}(Y)$ is an isomorphism for $k\leq m-2$, injective for $k=m-1$, and surjective for $k>m-1$. Define the fixed cohomology of $Y$ by $$H^{m-1}_{\textup{fixed}}(Y)=\textup{im}\left(H^{m-1}(X)\xrightarrow{\iota'{}^{*}} H^{m-1}(Y)\right)$$ and define the vanishing cohomology of $Y$ by $$H^{m-1}_{\textup{van}}(Y)=\textup{ker}\left(H^{m-1}(Y)\xrightarrow{\iota'_{!}}H^{m+1}(X)\right),$$ where  $\iota'_{!}$ is the Gysin morphism associated to $\iota'$. This group is called the vanishing cohomology group because, when viewing $Y$ as a general fiber of the Lefschetz fibration $f$, it is generated by the vanishing cycles of the fibration \cite[Lemma 2.26]{Voi2}. Moreover, since all vanishing cycles are conjugate under the monodromy action \cite[Corollary 3.24]{Voi2}, we have $H^{m-1}_{\textup{van}}(Y)=0$ if and only if all vanishing cycles vanish in rational cohomology.

There is an orthogonal decomposition $$H^{m-1}(Y)=H^{m-1}_{\textup{fixed}}(Y)\oplus H^{m-1}_{\textup{van}}(Y),$$ (see e.g., \cite[Lemma 2.26]{Voi2}).

Similarly, since $B$ is a complete intersection of two sections of an ample line bundle, the Lefschetz hyperplane theorem implies that the restriction map $H^k(X)\rightarrow H^{k}(B)$ is an isomorphism for $k\leq m-3$, injective for $k=m-2$ and surjective for $k>m-2$. We define the fixed cohomology of $B$ by $$H^{m-2}_{\textup{fixed}}(B)=\textup{im}\left(H^{m-2}(X)\xrightarrow{\iota^*} H^{m-2}(B)\right)$$ and the vanishing cohomology by $$H^{m-2}_{\textup{van}}(B)=\textup{ker}\left(H^{m-2}(B)\xrightarrow{\iota_{!}}H^{m+2}(X)\right),$$ where $\iota_{!}$ is the Gysin morphism associated to $\iota$. As in the case of $Y$, we have an orthogonal decomposition $$H^{m-2}(B)=H^{m-2}_{\textup{fixed}}(B)\oplus H^{m-2}_{\textup{van}}(B).$$

We will need the following lemmas in the next section:

\begin{lemma} \label{can} Let  $f:\tilde{X}\rightarrow \mb{P}^1$, $\pi:\tilde{X}\rightarrow X$, $\mathcal{L}$ and $E$ be as above. Then
      $$f^*\mathcal{O}_{\mb{P}^1}(1)\cong \pi^{*}\mathcal{L}\otimes \mathcal{O}(-E).$$
 \end{lemma}
 \begin{proof}
 Indeed, from the blowup morphism $\tilde{X}\xrightarrow{\pi}X$, we have the canonical bundle formula $K_{\tilde{X}}\cong \pi^{*}K_{X}\otimes \mathcal{O}(E)$, and if we regard $\tilde{X}$ as the hypersurface of $\mathbb{P}^1\times X$ defined by a section of the line bundle $f^*\mathcal{O}_{\mb{P}^1}(1)\otimes \pi^{*}\mathcal{L}$ (where by abuse of notation we denote the projection maps on $\mathbb{P}^1\times X$ by $f$ and $\pi$ as well), then by adjunction we have $$K_{\tilde{X}}\cong f^*\mathcal{O}_{\mb{P}^1}(-2)\otimes \pi^{*}K_{X} \otimes f^*\mathcal{O}_{\mb{P}^1}(1)\otimes \pi^{*}\mathcal{L}.$$ Thus by the comparison of the two expressions of $K_{\tilde{X}}$, the expression of $f^{*}\mathcal{O}_{\mb{P}^1}(1)$ follows. 
 \end{proof}
 
  The following lemma describes the action of cup product with $c_1\left(f^*\mathcal{O}_{\mb{P}^1}(1)\right)$ on the cohomology of $\tilde{X}$:
 
 \begin{lemma}\label{cup}
     Let  $f:\tilde{X}\rightarrow \mb{P}^1$, $\pi:\tilde{X}\rightarrow X$, $\mathcal{L}$, $E$ and $B$ be as above. If we write $H^k(\tilde{X})=H^k(X)\oplus H^{k-2}(B)$ for all $k$, then the morphism  $$H^k(\tilde{X})\xrightarrow{\cup c_1\left(f^*\mathcal{O}_{\mb{P}^1}(1)\right)} H^{k+2}(\tilde{X})$$ is given by $$(\alpha,\beta)\mapsto (\alpha\cup c_1(\mathcal{L})+\iota_{!}\beta,\ -\iota^*\alpha-\beta \cup c_1(\iota^*\mathcal{L})).$$
 \end{lemma}
 \begin{proof} Note that by Theorem \ref{Hodgeblowup}, the cohomology of $H^k(\tilde{X})$ fits into the following pushout diagram: 
\begin{center}
    \begin{tikzcd}
H^{k-4}(B) \arrow[r,"\iota_{!}"] \arrow[d,"\Psi"] & H^k(X) \arrow[d,"\pi^*"] \\
H^{k-2}(E) \arrow[r,"\tilde{\iota}_{!}"] & H^k(\tilde{X}), 
\end{tikzcd}
\end{center} where the horizontal maps are Gysin morphisms, and $\Psi$ is the cup product with the Euler class of the universal quotient bundle $Q$ of $\pi|_{E}^*N_{B/X}$ on $E\cong\mb{P}(N_{B/X})$. Since $N_{B/X}\cong \iota^*\mathcal{L}\oplus \iota^*\mathcal{L} $, by Grothendieck's relation of Chern classes, we have $$h^2+2c_1(\pi|_E^*\iota^*\mathcal{L})h+e(\pi|_{E}^*N_{B/X})=0,$$ where $h=c_1(\mathcal{O}_{E}(1))$. By applying the formula  $e(\mathcal{O}_{E}(-1))e(Q)=e(\pi|_{E}^*N_{B/X})$, we deduce that $e(Q)=2c_1(\pi|_E^*\iota^*\mathcal{L})+h.$ In other words, if we write $H^{k-2}(E)=H^{k-2}(B)\oplus \left(H^{k-4}(B)\cup h\right)$, then $$\Psi(\alpha)=(\alpha\cup (2c_1(\iota^*\mathcal{L})) ,\alpha\cup h ).$$ By this and the projection formula, for a $k$-form $\pi^*\alpha+\tilde{\iota}_{!}\pi|_{E}^*\beta\in H^k(\tilde{X})$ with $\alpha \in H^k(X)$ and $\beta \in H^{k-2}(B)$, we compute that 
\begin{equation*}
 \begin{split}
 (\pi^*\alpha+\tilde{\iota}_{!}\pi|_{E}^*\beta)\cup [-E]&=\tilde{\iota}_{!}\left(-\tilde{\iota}^*\pi^*\alpha+\pi|_{E}^*\beta\cup h\right)\\
&=\tilde{\iota}_{!}\left(-\pi|_{E}^*\iota^*\alpha+ \Psi(\beta)-2\pi|_{E}^*\left(\beta \cup c_1(\iota^*\mathcal{L})\right) \right)\\
&=\pi^*\iota_{!}\beta+ \tilde{\iota}_{!}\pi|_{E}^*\left(-\iota^*\alpha-2\beta \cup c_1(\iota^*\mathcal{L})\right).
 \end{split}   
\end{equation*}
In other words, if we write $H^k(\tilde{X})=H^k(X)\oplus H^{k-2}(B)$ for all $k$, then the morphism $$H^k(\tilde{X})\xrightarrow{\cup [-E]} H^{k+2}(\tilde{X})$$ is given by $$(\alpha,\beta)\mapsto (\iota_{!}\beta,-\iota^*\alpha-2\beta \cup c_1(\iota^*\mathcal{L})).$$ 

Therefore, using the isomorphism $$f^*\mathcal{O}_{\mb{P}^1}(1)\cong \pi^{*}\mathcal{L}\otimes \mathcal{O}(-E)$$ from Lemma \ref{can}, we comclude that the cup product  $$H^k(\tilde{X})\xrightarrow{\cup c_1\left(f^*\mathcal{O}_{\mb{P}^1}(1)\right)} H^{k+2}(\tilde{X})$$ acts via $$(\alpha,\beta)\mapsto (\alpha\cup c_1(\mathcal{L})+\iota_{!}\beta,\ -\iota^*\alpha-\beta \cup c_1(\iota^*\mathcal{L})).$$ 
\end{proof}

\section{Proof of the Main Theorem}

For $i=1,2$, let $X_i$ be a projective manifold of dimension $m_i$, $\{Y_{i,t}\}_{t\in \mathbb{P}^1}$ be a Lefschetz pencil on $X_i$ in the ample linear system $|\mathcal{L}_i|$. Let $$Y_i:=Y_{i,t} \textup{ (a general member)},\  B_i:= \textup{base locus of the pencil},\  \tilde{X}_i:=\textup{Bl}_{B_i}X_i,$$ and define $$E_i:= \textup{exceptional divisor},\  \iota_i:B_i\hookrightarrow X_i, \ \iota'_i:Y_i\hookrightarrow X_i, \ \tilde{\iota}_i: E_i\hookrightarrow \tilde{X}_i,\ \pi_i:\tilde{X}_i\rightarrow X_i,\  f_i:\tilde{X}_i\rightarrow\mathbb{P}^1$$ as in the previous section. We are only interested in the case that $\Delta(f_1)$ and $\Delta(f_2)$ are disjoint. In this case, it is easy to see that the fiber product $\tilde{X}_1\times_{\mathbb{P}^1}\tilde{X}_2$ is a smooth projective manifold. Now we are ready to prove Theorem \ref{fib}

\subsection{Mayer-Vietoris Argument}

\begin{proof}[Proof of Theorem \ref{fib}]

 We have a natural map $$\tilde{X}_1\times \tilde{X}_2\xrightarrow{\pi} X_1\times X_2$$ induced by the blowup morphisms $\tilde{X_i}\xrightarrow{\pi_i}X_i$ for $i=1,2$. Consider the image of the fiber product $\tilde{X}_1\times_{\mathbb{P}^1}\tilde{X}_2$ under this map, and denote it by $W$. Suppose that the Lefschetz pencils $\{Y_{i,t}\}_{t\in \mathbb{P}^1}$ are defined by the secitons $\sigma_{i,t}=\sigma_{i,0}+t\sigma_{i,\infty}\  ( t\in\mb{P}^1)$ of $|\mathcal{L}_i|$ on $X_i$. Then we can write $$W=V(\sigma_{1,0}\sigma_{2,\infty}-\sigma_{2,0}\sigma_{1,\infty})\subset X_1\times X_2.$$ Note that $W$ is singular precisely at the locus $$B_1\times B_2=V(\sigma_{1,0},\sigma_{1.\infty},\sigma_{2,0},\sigma_{2,\infty})\subset X_1\times X_2.$$ Moreover, locally at any point $p\in B_1\times B_2$, the germ $(W,p)$ is analytically isomorphic to the product of the germ of a threefold ordinary double point and $(B_1\times B_2,p)$. The fiber product $\tilde{X}_1\times_{\mathbb{P}^1}\tilde{X}_2$ can be regarded as a small resolution of $W$ with exceptional set $E$. Here, $E$ is a $\mathbb{P}^1$-bundle over the base $Z:=B_1\times B_2$ isomorphic to $$\mathbb{P}(p_1^*\iota_1^*\mathcal{L}_1\oplus p_2^*\iota_2^*\mathcal{L}_2),$$ where $p_i:B_1\times B_2\rightarrow B_i$ denote the natural projection map. 
 
 Thus, we have the following discriminant square 
\begin{center}
    \begin{tikzcd}
E  \arrow[r,"\tilde{\iota}"] \arrow[d, "\pi|_{E}"] & \tilde{X}_1\times_{\mathbb{P}^1}\tilde{X}_2 \arrow[d, "\pi|_{\tilde{X}_1\times_{\mathbb{P}^1}\tilde{X}_2}"] \\
Z \arrow[r,"\iota"]                               & W,                                                             
\end{tikzcd}
\end{center} which naturally embeds into the larger discriminant square associated to the birational map $\pi$: 

\begin{center}
    \begin{tikzcd}
E':= E_1\times \tilde{X}_2 \cup \tilde{X}_1\times E_2  \arrow[r,"\tilde{\iota}'"] \arrow[d,"\pi|_{E'}"] & \tilde{X}_1\times \tilde{X}_2 \arrow[d, "\pi"] \\
Z':=B_1\times X_2\cup X_1\times B_2 \arrow[r,"\iota'"]                               & X_1\times X_2.                                                              
\end{tikzcd}
\end{center}
Note that the two components of $E'$ intersect along $E_1\times E_2$, and $E$ can be identified with the diagonal subvariety inside this intersection. 

Then, by the Mayer-Vietoris sequence for the discriminant square (cf. \cite[Corollary-Definition 5.37]{PS}), we obtain a commutative diagram of long exact sequences of mixed Hodge structures: 

\begin{center}
    \begin{tikzcd}
... \arrow[r] & H^k(X_1\times X_2) \arrow[r] \arrow[d] & H^k(\tilde{X}_1\times \tilde{X}_2)\oplus H^k(Z')  \arrow[d] \arrow[r]                    & H^k(E') \arrow[r] \arrow[d] & ... \\
... \arrow[r] & H^k(W) \arrow[r]                       & H^k(\tilde{X}_1\times_{\mathbb{P}^1}\tilde{X}_2)\oplus H^k(Z)  \arrow[r] & H^k(E) \arrow[r]            & ...
\end{tikzcd}
\end{center}

Since $E'$ is compact and $X_1\times X_2$ is smooth, the connecting homomorphism $H^{k-1}(E')\rightarrow H^k(X_1\times X_2)$ vanished for all $k$. Thus, the top row in the diagram is a short exact sequence for each $k$. We now analyze the vertical maps in the diagram above:
\begin{enumerate}
    \item \textbf{Map $H^k(X_1\times X_2)\rightarrow H^k(W)$:} 
    
    The variety $W$ is a section of the ample line bundle $q_1^{*}\mathcal{L}_1\otimes q_2^*\mathcal{L}_2$ on $X_1\times X_2$, where $q_i:X_1\times X_2\rightarrow X_i$ is the projection. By the Lefschetz hyperplane theorem (cf. \cite[Theorem C.15]{PS}), the restriction map $$H^k(X_1\times X_2)\rightarrow H^k(W)$$ is an isomorphism for $k\leq m_1+m_2-2$, and injective for $k=m_1+m_2-1$. In particular, the mixed Hodge structure on $H^k(W)$ is pure for $k\leq m_1+m_2-2$. Consequently, for these values of $k$, the bottom horizontal sequence is injective on the left, and for $k\leq m_1+m_2-3$, it is surjective on the right.

    \item \textbf{Maps $H^k(E')\rightarrow H^k(E)$ and $H^k(Z')\rightarrow H^k(Z)$:} 

We can analyze $H^k(E')$ and $H^k(Z')$ using the Mayer-Vietoris sequence associated to their two irreducible components. By comparing the two Mayer-Vietoris sequences, one can obtain the identification $$H^k(E')\cong H^k(Z')\oplus H^{k-2}(B_1\times X_2)\oplus H^{k-2}(X_1\times B_2)\oplus H^{k-4}(B_1\times B_2).$$ Write $H^k(E)=H^k(Z)\oplus H^{k-2}(Z)$. Then the map $H^k(E')\rightarrow H^k(E)$ can be described as follows: 
\begin{enumerate}
    \item The summand $H^k(Z')$ maps to $H^k(Z)$ via the natural restriction map $H^k(Z')\rightarrow H^k(Z)$, and maps to the summand $H^{k-2}(Z)$ via the zero map. 
    \item $$H^{k-2}(B_1\times X_2)\oplus H^{k-2}(X_1\times B_2)$$ maps to $H^{k}(Z)$ via the zero map, and maps to $H^{k-2}(B_1\times B_2)=H^{k-2}(Z)$ via restriction.
    \item The summand  $H^{k-4}(B_1\times B_2)$ maps as follows: any $\alpha \in H^{k-4}(B_1\times B_2)$ is sent to $$-\pi|_E^*\alpha \cup c_1(\mathcal{O}_E(1))^2\in H^k(E),$$ which, by Grothendieck's relation of Chern classes, yields the map $$\alpha\mapsto \left(\alpha\cup c_1(q_1^{*}\iota_1^*\mathcal{L}_1)\cup c_1(q_2^{*}\iota_2^*\mathcal{L}_2),\alpha \cup \left(c_1(q_1^{*}\iota_1^*\mathcal{L}_1)+ c_1(q_2^{*}\iota_2^*\mathcal{L}_2)\right) \right)$$ into $H^k(Z)\oplus H^{k-2}(Z)$
\end{enumerate}

It is not hard to see that the image of $H^{k-4}(B_1\times B_2)$ is contained in the image of the other three direct summands of $H^k(E')$ for all $k$. Hence, the dimension of $$\textup{ker}(H^k(E')\rightarrow H^k(E))$$ is $$ \textup{dim}_{\C}\textup{ker}(H^k(Z')\rightarrow H^k(Z)) + \textup{dim}_{\C}H^{k-4}(B_1\times B_2)$$$$ +\textup{dim}_{\C}\textup{ker}( H^{k-2}(B_1\times X_2)\oplus H^{k-2}(X_1\times B_2)\rightarrow H^{k-2}(B_1\times B_2))$$ for all $k$. 

Moreover, we note that $$\textup{coker}(H^k(E')\rightarrow H^k(E))\cong\textup{coker}(H^k(Z')\rightarrow H^k(Z))$$ for all $k\neq m_1+m_2-4$, and  \begin{equation*}
    \begin{split}
        \textup{coker}(H^{m_1+m_2-4}(E')\rightarrow H^{m_1+m_2-4}(E))\cong&\textup{coker}(H^{m_1+m_2-4}(Z')\rightarrow H^{m_1+m_2-4}(Z))\\ &\oplus H^{m_1-2}_{\textup{van}}(B_1)\otimes H^{m_2-2}_{\textup{van}}(B_2).
    \end{split}
\end{equation*}

\end{enumerate}

 Now, returning to the cohomology of the fiber product $H^k(\tilde{X}_1\times_{\mathbb{P}^1}\tilde{X}_2)$, we prove the following dimension estimates:
 
\begin{lemma} \label{dimbound}
The degree $k$ part of the kernel of the quotient map $$H^{\bullet}(\tilde{X}_1\times \tilde{X}_2)=H^{\bullet}(\tilde{X}_1)\otimes_{\C} H^{\bullet}(\tilde{X}_2)\rightarrow H^{\bullet}(\tilde{X}_1)\otimes_{H^{\bullet}(\mb{P}^1)}H^{\bullet}(\tilde{X}_2)$$ has dimension at least $$\textup{dim}_{\C}H^{k-4}(B_1\times B_2)$$$$ +\textup{dim}_{\C}\textup{ker}( H^{k-2}(B_1\times X_2)\oplus H^{k-2}(X_1\times B_2)\rightarrow H^{k-2}(B_1\times B_2))$$ for all $k\leq m_1+m_2-1$.

\end{lemma}
 
 Assuming the lemma, then by the Snake Lemma, the canonical map $$H^{\bullet}(\tilde{X}_1)\otimes_{H^{\bullet}(\mb{P}^1)}H^{\bullet}(\tilde{X}_2)\rightarrow H^{\bullet}(\tilde{X}_1\times_{\mathbb{P}^1}\tilde{X}_2)$$ is injective in degree $k\leq m_1+m_2-1$, and subjective in degree $k\leq m_1+m_2-3$. Moreover, the dimension bound of the kernel in Lemma \ref{dimbound} is an equality.

\begin{proof}[Proof of Lemma \ref{dimbound}]

 We begin by writing the cohomology of the product $\tilde{X}_1 \times \tilde{X}_2$ as $$H^k(\tilde{X}_1\times \tilde{X}_2)=H^k(X_1\times X_2)\oplus H^{k-2}(B_1\times X_2)\oplus H^{k-2}(X_1\times B_2)\oplus H^{k-4}(B_1\times B_2).$$ By applying Lemma \ref{cup}, it is not hard to see that the degree $k$ part of the kernel of the quotient map $$H^{\bullet}(\tilde{X}_1\times \tilde{X}_2)\rightarrow H^{\bullet}(\tilde{X}_1)\otimes_{H^{\bullet}(\mb{P}^1)}H^{\bullet}(\tilde{X}_2)$$ is generated by  $4$-tuples $(M_1,M_2,M_3,M_4)$, where $$M_1=(\alpha\cup c_1(\mathcal{L}_1)+\iota_{1!}\beta)\cup \gamma -\alpha\cup (\gamma\cup c_1(\mathcal{L}_2)+\iota_{2!}\delta),$$ $$M_2=(-\iota_1^*\alpha-\beta \cup c_1(\iota_1^*\mathcal{L}_1))\cup \gamma - \beta \cup (\gamma\cup c_1(\mathcal{L}_2)+\iota_{2!}\delta),$$ $$M_3= (\alpha\cup c_1(\mathcal{L}_1)+\iota_{1!}\beta)\cup \delta-\alpha \cup (-\iota_2^*\gamma-\delta \cup c_1(\iota_2^*\mathcal{L}_2)),$$ $$M_4=(-\iota_1^*\alpha-\beta \cup c_1(\iota_1^*\mathcal{L}_1))\cup \delta -\beta \cup (-\iota_2^*\gamma-\delta \cup c_1(\iota_2^*\mathcal{L}_2)),$$ where $\alpha\in H^l(X_1),$ $\beta\in H^{l-2}(B_1)$, $\gamma\in H^{k-l-2}(X_2)$. $\delta\in H^{k-l-4}(B_2)$ and $0\leq l\leq k$. 
 
  Set $\beta=\delta=0$. Then $$(M_1,M_2,M_3,M_4)=(\alpha\cup c_1(\mathcal{L}_1)\cup\gamma-\alpha\cup\gamma\cup c_1(\mathcal{L}_2),-\iota_1^*\alpha\cup \gamma, \alpha\cup\iota_2^*\gamma,0).$$ Projecting onto the $(M_2, M_3)$ summand, the image spans a subspace whose dimension is at least $$ \textup{dim}_{\C}\textup{ker}( H^{k-2}(B_1\times X_2)\oplus H^{k-2}(X_1\times B_2)\rightarrow H^{k-2}(B_1\times B_2)).$$ 
  
  Now we set either $\beta=\gamma=0$ or $\alpha=\delta=0$. Then projecting to the $M_4$ summand, we see that these elements generate $$\textup{im}\left(H^{k-4}(X_1\times B_2)\oplus H^{k-4}(B_1\times X_2)\rightarrow H^{k-4}(B_1\times B_2)\right).$$  This image has dimension at least  $$\textup{dim}_{\C}H^{k-4}(B_1\times B_2)$$ for $k\neq m_1+m_2$, as the corresponding maps are surjective in those degrees. 
  Combining the two dimension bounds completes the proof of the lemma. 
\end{proof}

The surjectivity in $k=m_1+m_2-2$ doesn't follow from this picture. Instead, we will use the Leray spectral sequence to compute the dimension of $H^{m_1+m_2-2}(\tilde{X}_1\times_{\mb{P}^1}\tilde{X}_2)$ and show that it is equal to the dimension of the degree $m_1+m_2-2$ part of $$H^{\bullet}(\tilde{X}_1)\otimes_{H^{\bullet}(\mb{P}^1)}H^{\bullet}(\tilde{X}_1).$$

We will prove the following lemma in the next subsection, where $h^k$, $h^k_{\textup{van}}$ denote the $\Q$-dimension of the cohomology $H^k$, $H^k_{\textup{van}}$ respectively:

\begin{lemma}\label{dim}
    Let $d_i$ be the cardinality of $\Delta(f_i)$. Then \begin{enumerate}
        \item $h^{m_i}(\tilde{X}_i)=d_i+h^{m_i}(Y_i)+h^{m_i-2}(Y_i)-2h^{m_i-1}_{\textup{van}}(Y_i),$
        \item The following formula holds: \begin{equation*}
            \begin{split}
                h^{m_1+m_2-2}(\tilde{X}_1\times_{\mb{P}^1}\tilde{X}_2)=& h^{m_1+m_2-2}(Y_1\times Y_2)-h_{\textup{van}}^{m_1-1}(Y_1)h^{m_2-1}(Y_2)\\&-h^{m_1-1}(Y_1)h^{m_2-1}_{\textup{van}}(Y_2)+h^{m_1-1}_{\textup{van}}(Y_1)h^{m_2-1}_{\textup{van}}(Y_2)\\
                &+h^{m_1+m_2-4}(Y_1\times Y_2)-h_{\textup{van}}^{m_1-1}(Y_1)h^{m_2-3}(Y_2)\\ &-h^{m_1-3}(Y_1)h^{m_2-1}_{\textup{van}}(Y_2)+d_1 h^{m_2-2}(Y_2)+d_2 h^{m_1-1}(Y_1)\\&-2\left(h^{m_1-1}_{\textup{van}}(Y_1)h^{m_2-2}(Y_2)+h^{m_1-2}(Y_1)h^{m_2-1}_{\textup{van}}(Y_2)\right)
            \end{split}
        \end{equation*}
    \end{enumerate}
\end{lemma}

Assuming the lemma, note also that the degree $m_1+m_2-2$ part of  $$H^{\bullet}(\tilde{X}_1)\otimes_{H^{\bullet}(\mb{P}^1)}H^{\bullet}(\tilde{X}_1)$$ has dimension \begin{equation}\label{dim'}
    \begin{split}
        h^{m_1+m_2-2}(\tilde{X}_1\times \tilde{X}_2)-h^{m_1+m_2-6}(B_1\times B_2)
        -h^{m_1+m_2-4}(B_1\times X_2)\\
        -h^{m_1+m_2-4}(X_1\times B_2)+h^{m_1+m_2-4}(B_1\times B_2)-h^{m_1-2}_{\textup{van}}(B_1)h^{m_2-2}_{\textup{van}}(B_2).
    \end{split}
\end{equation}
Thus, by canceling $d_i$ of Lemma \ref{dim} (2) using Lemma \ref{dim} (1), it is pretty elementary to verify that $h^{m_1+m_2-2}(\tilde{X}_1\times_{\mb{P}^1}\tilde{X}_2)$ is equal to the formula (\ref{dim'}) by the Lefschetz hyperplane theorem, the K\"unneth formula and the Poincar\'e duality. This proves that the canonical map $$H^{\bullet}(\tilde{X}_1)\otimes_{H^{\bullet}(\mb{P}^1)}H^{\bullet}(\tilde{X}_2)\rightarrow H^{\bullet}(\tilde{X}_1\times_{\mathbb{P}^1}\tilde{X}_2)$$ induce an isomorphism in degree $k=m_1+m_2-2$, and the proof is complete.
\end{proof}

Corollary \ref{hformula1} is an immediate consequence of the computation in the proof of Lemma \ref{dimbound}:

\begin{proof}[Proof of Corollary \ref{hformula1}]
By Theorem \ref{fib}, it suffices to compute the Hodge numbers of the Hodge structures on the corresponding degree of $$H^{\bullet}(\tilde{X}_1)\otimes_{H^{\bullet}(\mb{P}^1)}H^{\bullet}(\tilde{X}_2).$$ The computation in the proof of Lemma \ref{dimbound} and the discussion after Lemma \ref{dimbound} not only show the dimension equality but also prove that there is an isomorphism of Hodge structures between the degree $k$ part of $$H^{\bullet}(\tilde{X}_1)\otimes_{H^{\bullet}(\mb{P}^1)}H^{\bullet}(\tilde{X}_2)\rightarrow H^{\bullet}(\tilde{X}_1\times_{\mathbb{P}^1}\tilde{X}_2)$$ and $$H^{k-4}(B_1\times B_2)\oplus\textup{ker}( H^{k-2}(B_1\times X_2)\oplus H^{k-2}(X_1\times B_2)\rightarrow H^{k-2}(B_1\times B_2))$$ for all $k\leq m_1+m_2-1$ (where the natural Tate twists are ignored for simplicity). The formula in Corollary \ref{hformula1} follows.
    
\end{proof}

\subsection{Dimension Counting via Leray Spectral Sequence}

\begin{proof}[Proof of Lemma \ref{dim}]
    We will apply the method discussed in \cite[Section 4.5.4]{PS}. We first consider the Lefschetz fibration $\tilde{X}_i\xrightarrow{f_i} \mb{P}^1$ and prove the formula (1). In this case, the Leray spectral sequence for $f_i$ degenerates at $E_2$ (cf. \cite[Theorem 4.24]{PS}). Moreover, the Lefschetz fibration satisfies the local invariant cycle property (\cite[Corollary C.21]{PS}): the adjunction homomorphism $$a_k: R^{k}f_{i*}\underline{\Z}_{X_i}\rightarrow j_{i*}j_i^{*}R^{k}f_{i*}\underline{\Z}_{X_i}$$ is surjective for all $k$, where $j_i:U_i:=\mb{P}^1-\Delta(f_i)\hookrightarrow \mb{P}^1$ is the inclusion map. 
    
    We need to divide the computation into two cases. 
    
    \textbf{Case I:}
    The first case (the generic case) is that $H^{m_i-1}_{\textup{van}}(Y_i)$ is non-zero, or equivalently, all of the vanishing cycles $\delta_x$ for $x\in \Delta(f_i)$ are non-zero in $H^{m_i-1}(Y_i)$ since the vanishing cycles are conjugate under monodromy. In this case, $R^{k}f_{i*}\underline{\Z}_{X_i}$ is locally free for $k\neq m_i-1$ and the adjunction homomorphism $a_{m_i-1}$ is an isomorphism (\cite[Lemma C.13]{PS}). Therefore, the Leray spectral sequence for $f_i$ reads $$E^{0,m_i}_2(f_i)=H^{m_i}(Y_i,\C),$$ $$E^{1,m_i-1}_2(f_i)=H^1(\mb{P}^1,j_{i*}j_i^{*}R^{k}f_{i*}\underline{\C}_{X_i}),$$ $$E^{2,m_i-2}_2(f_i)=H^{m_i-2}(Y_i,\C).$$ Let $\mathcal{L}$ be the local system $j_i^{*}R^{m_i-1}f_{i*}\underline{\C}_{X_i}$ on $U_i$. We can calculate $E^{1,m_i-1}_2(f_i)$ by the exact sequence $$0\rightarrow j_{i!}\mathcal{L}\rightarrow j_{i*}\mathcal{L}\rightarrow \bigoplus_{x\in\Delta(f_i)}\mathcal{L}^{T}_x\rightarrow 0,$$ where $\Delta^{*}_x$ is a small punctured disk centered at $x$, $\mathcal{L}_x$ is the restriction of $\mathcal{L}$ on $\Delta^{*}_x$, and $\mathcal{L}_x^T$ is the subspace of invariants in $\mathcal{L}_x$ under the local monodromy $T$. By taking the long exact sequence of the cohomology, we have $$0\rightarrow H^{0}(U_i,\mathcal{L})\rightarrow \bigoplus_{x\in\Delta(f_i)}\mathcal{L}^{T}_x \rightarrow H^1_c(U_i,\mathcal{L})\rightarrow H^1(\mb{P}^1,j_{i*}\mathcal{L})\rightarrow 0.$$ By the theorem of the fixed part (\cite[Theorem 4.23]{PS}), $$H^0(U_i,\mathcal{L})=\textup{im}(H^{m_i-1}(\tilde{X}_i)\rightarrow H^{m_i-1}(Y_i))=\textup{im}(H^{m_i-1}(X_i)\rightarrow H^{m_i-1}(Y_i)).$$ The last equality follows from Lemma \ref{gysin} since the Gysin morphism $H^{m_i-3}(B_i)\rightarrow H^{m_i-1}(Y_i)$ can be identified via the Lefschetz hyperplane theorem with the composition $H^{m_i-3}(X_i)\rightarrow H^{m_i-1}(X_i)\rightarrow H^{m_i-1}(Y_i)$, where the first map is the Lefschetz operator given by cup product with the hyperplane class, and the second is the restriction map. Hence, we see that $$h^{0}(U,\mathcal{L})=h^{m_i-1}(Y_i)-h^{m_i-1}_{\textup{van}}(Y_i).$$ Also, since $\mathcal{L}^T_{x}$ is generated by the class in $H^{m_i-1}(Y_i)$ which is orthogonal to the vanishing cycle $\delta_x$, $\mathcal{L}^T_{x}$ has dimension $h^{m_i-1}(Y_i)-1$ for all $x\in \Delta(f_i)$. On the other hand, $h^1_c(U_i,\mathcal{L})=h^1(U_i,\mathcal{L}^{\vee})$ by the Poincar\"e duality and we have $$H^k(U_i,\mathcal{L}^{\vee})=\textup{Ext}^k_{\Z[\pi_1(U_i)]}(\Z,\mathcal{L}^{\vee})$$ for all $k$. Since $U_i$ is homotopic to $\bigvee_{i=1}^{d_i-1}S^1$, $\Z[\pi_1(U_i)]$ is the group ring of a free group of $d_i-1$ generators. There is a natural two-term free resolution of $\Z$ as a $\Z[\pi_1(U_i)]$ module (cf. \cite[Corollary 6.2.7]{weibel}): $$0\rightarrow \Z[\pi_1(U_i)]^{d_i-1}\rightarrow \Z[\pi_1(U_i)] \rightarrow \Z\rightarrow 0.$$ Therefore, by the isomorphism $H^{0}(U_i,\mathcal{L}^{\vee})\cong H^{0}(U_i,\mathcal{L})$, we have $$h^1(U_i,\mathcal{L}^{\vee})-h^0(U_i,\mathcal{L})=(d_i-2)\textup{rank}(\mathcal{L}).$$ As a consequence, \begin{equation*}
        \begin{split}
            h^{1}(\mb{P}^1,j_{i*}\mathcal{L})&=h^1(U_i,\mathcal{L}^{\vee})+h^0(U_i,\mathcal{L})-\textup{dim}\left(\bigoplus_{x\in\Delta(f_i)}\mathcal{L}_x^T\right)\\
            &= 2h^{0}(U_i,\mathcal{L})+(d_i-2)\textup{rank}(\mathcal{L})-d_i(\textup{rank}\mathcal{L}-1)\\
            &=d_i-2h^{m_i-1}_{\textup{van}}(Y_i).
        \end{split}
    \end{equation*}
     By summing the dimension of $E^{0,m_i}(f_i)$, $E^{1.m_i-1}(f_i)$ and $E^{2,m_i-2}(f_i)$, we get the formula in Lemma \ref{dim} (1).
     
 \textbf{Case II:}
   The second case we need to consider is the special case where $H^{m_i-1}_{\textup{van}}(Y_i)$ is zero, i.e., the restriction map $H^{m_i-1}(X_i) \rightarrow H^{m_i-1}(Y_i)$ is an isomorphism. In this case, $R^{k}f_{i*}\underline{\C}_{X_i}$ is locally constant for $k\neq m_i$, and for $k=m_i$, we have the short exact sequence $$0\rightarrow \bigoplus_{x\in \Delta(f_i)} H^{m_i-1}(F_x)\rightarrow R^{m_i}f_{i*}\underline{\C}_{X_i}\rightarrow j_{i*}j_{i}^{*}R^{m_i}f_{i*}\underline{\C}_{X_i}\rightarrow 0,$$ where $H^{m_i-1}(F_x)$ is the skyscraper sheaf (supported on $x$) of the cohomology of the Milnor fiber of $f_i$ around $x$. Also, $j_{i*}j_{i}^{*}R^{m_i}f_{i*}\underline{\C}_{X_i}$ is a locally constant sheaf with stalks isomorphic to $H^{m_i}(Y_i)$. Therefore, the Leray spectral sequence for $f_i$ reads  $$E^{0,m_i}_2(f_i)=H^{0}(\mb{P}^1,R^{m_i}f_{i*}\underline{\C}_{X_i}),$$ $$E^{1,m_i-1}_2(f_i)=0,$$ $$E^{2,m_i-2}_2(f_i)=H^{m_i-2}(Y_i,\C),$$ and we have $$h^{0}(\mb{P}^1,R^{m_i}f_{i*}\underline{\C}_{X_i})=d_i+h^{m_i}(Y_i,\C).$$ The formula in Lemma \ref{dim} (1) follows.

   Now we turn to the fibration $f:\tilde{X}_1\times_{\mathbb{P}^1} \tilde{X}_2\rightarrow \mb{P}^1$ and prove the formula in Lemma \ref{dim} (2). Write $\Delta(f)=\Delta(f_1)\cup\Delta(f_2)$. We want to carry out the same computation as in the Lefschetz fibration case for $f$. Since $\Delta(f_1)$ and $\Delta(f_2)$ are disjoint and the local monodromy action on $H^k(Y_1\times Y_2)=\bigoplus_{i+j=k}H^i(Y_1)\otimes H^j(Y_2)$ is induced by product group action of local monodromy 
 group action on $H^i(Y_1)$ and $H^j(Y_2)$, we still have the local invariant cycle property for $f$: the adjunction homomorphism $$a_k: R^{k}f_{*}\underline{\Z}_{\tilde{X}_1\times_{\mb{P}^1}\tilde{X}_2}\rightarrow j_{*}j^{*}R^{k}f_{*}\underline{\Z}_{\tilde{X}_1\times_{\mb{P}^1}\tilde{X}_2}$$ is surjective for all $k$, where $j:U:=\mb{P}^1-\Delta(f)\hookrightarrow \mb{P}^1$ is the inclusion map. 
   
We will now analyze the behavior of the Leray spectral sequence in three cases, depending on the vanishing cohomology of the fibers $Y_1$ and $Y_2$:

\textbf{Case I$'$}: Both $Y_1$ and $Y_2$ have non-trivial vanishing cohomology,
    
    \textbf{Case II$'$}: Exactly one on then does,
   
    \textbf{Case III$'$}: Neither has non-trivial vanishing cohomology.

Let's start with the generic case (Case I$'$).

\textbf{Case I$'$:} $H^{m_i}_{\textup{van}}(Y_i)\neq 0$ for $i=1,2$. 

For $x\in \Delta(f_1)$, let $\Delta_x$ be a small punctured disc centered at $x$ and $B_x$ be a small ball around the singularity at $f_1^{-1}(x)$. Let $B_x\cap f_1^{-1}(\Delta_x-\{x\})\xrightarrow{f_1}\Delta_x-\{x\}$ denote the Milnor fibration around the singularity at $f_1^{-1}(x)$ and $F_x$ the Milnor fiber. By the K\"unneth formula and the excision theorem, for all $t\in \Delta_x-\{x\}$, we have \begin{equation*}
       \begin{split}
           H^{k+1}(f^{-1}\Delta_x,f^{-1}(t))&=H^{k+1}(f_1^{-1}\Delta_x\times Y_2, Y_1\times Y_2)\\
           &\cong \bigoplus_{i+j=k+1}H^i(f_1^{-1}\Delta_x,Y_1)\otimes H^j(Y_2)\\
           &\cong \bigoplus_{i+j=k+1}H^i(B_x,F_x)\otimes H^j(Y_2)\\
           &\cong H^{m_1-1}(F_x)\otimes H^{k+1-m_1}(Y_2).   
       \end{split}
 \end{equation*} 
   Therefore, the long exact sequence of the cohomology for the pair $(f^{-1}\Delta_x,f^{-1}(t))$ gives $$0\rightarrow H^k(f^{-1}\Delta_x)\rightarrow H^k(Y_1\times Y_2)\rightarrow H^{m_1-1}(F_x)\otimes H^{k+1-m_1}(Y_2)\rightarrow 0.$$ The surjectivity on the right (and the injectivity on the left) follows from the assumption that the vanishing cycles for $f_1$ survive globally. A similar short exact sequence holds for $x\in \Delta(f_2)$. In particular, the adjunction homomorphism  $$a_k: R^{k}f_{*}\underline{\C}_{\tilde{X}_1\times_{\mb{P}^1}\tilde{X}_2}\rightarrow j_{*}j^{*}R^{k}f_{*}\underline{\C}_{\tilde{X}_1\times_{\mb{P}^1}\tilde{X}_2}$$ is an isomorphism for all $k$. For all $k$, let $\mathcal{L}^k$ be the local system $j^*R^{k}f_{*}\underline{\C}_{\tilde{X}_1\times_{\mb{P}^1}\tilde{X}_2}$ on $U$. We have the short exact sequence $$0\rightarrow j_{!}\mathcal{L}^k\rightarrow j_{*}\mathcal{L}^k\rightarrow \bigoplus_{x\in\Delta(f)}(\mathcal{L}^k_x)^{T}\rightarrow 0,$$ where $\Delta^{*}_x$ is a small punctured disk centered at $x$, $\mathcal{L}^k_x$ is the restriction of $\mathcal{L}^k$ on $\Delta^{*}_x$, and $(\mathcal{L}^k_x)^T$ is the subspace of invariants in $\mathcal{L}^k_x$ under the local monodromy $T$. For each $k$, $$H^0(\mb{P}^1,j_{*}\mathcal{L}^k)=H^k(Y_1\times Y_2)^{\pi_1(U)}=\bigoplus_{i+j=k}H^i(Y_1)^{\pi_1(U_1)}\otimes H^j(Y_2)^{\pi_1(U_2)},$$ and for $i=1,2$, $$H^l(Y_i)^{\pi_1(U_i)}=\left\{ \begin{array}{cc}
      H^{l}(Y_i),  & l\neq m_i-1 \\
       \textup{im}\left(H^l(X_i)\rightarrow H^l(Y_i)\right), & l=m_i-1.
   \end{array}\right.$$ Also, $H^2(\mb{P}^1,j_*\mathcal{L}^k)=H^2(\mb{P}^1,j_!\mathcal{L}^k)=H^2_c(U_i,\mathcal{L}^k)$ is dual to $H^0(U,(\mathcal{L}^k)^{\vee})=H^0(\mb{P}^1,j_{*}\mathcal{L}^k)$. So we have $$h^2(\mb{P}^1,j_*\mathcal{L}^k)=h^0(\mb{P}^1,j_*\mathcal{L}^k).$$
    On the other hand, we have the exact sequence $$0\rightarrow H^{0}(U,\mathcal{L}^k)\rightarrow \bigoplus_{x\in\Delta(f)}(\mathcal{L}^k_x)^{T} \rightarrow H^1_c(U,\mathcal{L}^k)\rightarrow H^1(\mb{P}^1,j_{*}\mathcal{L}^k)\rightarrow 0.$$ Note that $$\textup{dim}(\mathcal{L}^k_x)^T=\left\{\begin{array}{cc}
    h^{k}(Y_1\times Y_2)-h^{k+1-m_1}(Y_2), & x\in \Delta(f_1)\\
    h^{k}(Y_1\times Y_2)-h^{k+1-m_2}(Y_1), & x\in \Delta(f_2)
    \end{array}\right.,$$ and similar to the Lefschetz fibration case, we have $$h^{1}_c(U,\mathcal{L}^k)-h^{0}(U,\mathcal{L}^k)=(d_1+d_2-2)\textup{rank}\mathcal{L}^k.$$ As a result, \begin{equation*}
        \begin{split}
            h^1(\mathbb{P}^1,j_{*}\mathcal{L}^k)&=h^1_c(U,\mathcal{L}^k)+h^0(U,\mathcal{L}^k)-\textup{dim}\left(\bigoplus_{x\in\Delta(f)}(\mathcal{L}^k_x)^T\right)\\
            &=2h^{0}(U,\mathcal{L}^k)+(d_1+d_2-2)\textup{rank}\mathcal{L}^k\\ &\ \ \ -d_1(\textup{rank}\mathcal{L}^k-h^{k+1-m_1}(Y_2))-d_2(\textup{rank}\mathcal{L}^k-h^{k+1-m_2}(Y_1))\\
            &=2h^0(U,\mathcal{L}^k)+d_1h^{k+1-m_1}(Y_2)+d_2h^{k+1-m_2}(Y_1)-2h^k(Y_1\times Y_2)
        \end{split}
    \end{equation*}
    
    It follows that \begin{equation*}
        \begin{split}
            \textup{dim}_{\C}E^{0,m_1+m_2-2}_2(f)&=h^{0}(\mb{P}^1,j_{*}\mathcal{L}^{m_1+m_2-2})\\
            &=h^{m_1+m_2-2}(Y_1\times Y_2)-h_{\textup{van}}^{m_1-1}(Y_1)h^{m_2-1}(Y_2)\\&\ \ \ -h^{m_1-1}(Y_1)h^{m_2-1}_{\textup{van}}(Y_2)+h^{m_1-1}_{\textup{van}}(Y_1)h^{m_2-1}_{\textup{van}}(Y_2),
        \end{split}
    \end{equation*}
\begin{equation*}
        \begin{split}
            \textup{dim}_{\C}E^{1,m_1+m_2-3}_2(f)&=h^{1}(\mb{P}^1,j_{*}\mathcal{L}^{m_1+m_2-3})\\
            &=d_1h^{m_2-2}(Y_2)+d_2h^{m_1-2}(Y_1)\\&\ \ \ -2\left(h^{m_1-1}_{\textup{van}}(Y_1)h^{m_2-2}(Y_2)+h^{m_1-2}(Y_1)h^{m_2-1}_{\textup{van}}(Y_2)\right),
        \end{split}
    \end{equation*}
    
and \begin{equation*}
        \begin{split}
            \textup{dim}_{\C}E^{2,m_1+m_2-4}_2(f)&=h^{0}(\mb{P}^1,j_{*}\mathcal{L}^{m_1+m_2-4})\\
            &=h^{m_1+m_2-4}(Y_1\times Y_2)-h_{\textup{van}}^{m_1-1}(Y_1)h^{m_2-3}(Y_2)\\&\ \ \ -h^{m_1-3}(Y_1)h^{m_2-1}_{\textup{van}}(Y_2).
        \end{split}
    \end{equation*}
We get the formula in Lemma \ref{dim} (2) by summing the three dimensions. 

\textbf{Case II$'$:} $H^{m_1-1}_{\textup{van}}(Y_1)= 0$ and $H^{m_2-1}_{\textup{van}}(Y_2)\neq 0$. 

In this case, for $x\in\Delta(f_2)$, the same local computation of $R^kf_{*}\underline{\C}_{\tilde{X}_1\times_{\mb{P}^1}\tilde{X}_2}$ as in the previous case holds. For $x\in\Delta(f_1)$, since the cohomology class of the  vanishing cycle $\delta_x$ is zero, we have instead the short exact sequence $$0\rightarrow H^{m_1-1}(F_x)\otimes H^{k-m_1}(Y_2)\rightarrow H^k(f^{-1}\Delta_x)\rightarrow H^k(Y_1\times Y_2)\rightarrow 0.$$ Thus we have a short exact sequence $$0\rightarrow \bigoplus_{x\in\Delta(f_1)}H^{m_1-1}(F_x)\otimes H^{k-m_1}(Y_2)\rightarrow R^kf_{*}\underline{\C}_{\tilde{X}_1\times_{\mb{P}^1}\tilde{X}_2}\rightarrow j^*j_{*}R^kf_{*}\underline{\C}_{\tilde{X}_1\times_{\mb{P}^1}\tilde{X}_2}\rightarrow 0,$$ where $H^{m_1-1}(F_x)\otimes H^{k-m_1}(Y_2)$ is regarded as a skyscraper sheaf supported on $x$. Also, with the notation in the previous case, we have instead $$\textup{dim}(\mathcal{L}^k_x)^T=\left\{\begin{array}{cc}
    h^{k}(Y_1\times Y_2), & x\in \Delta(f_1)\\
    h^{k}(Y_1\times Y_2)-h^{k+1-m_2}(Y_1), & x\in \Delta(f_2)
    \end{array}\right.,$$
As a result, 
\begin{equation*}
        \begin{split}
            \textup{dim}_{\C}E^{0,m_1+m_2-2}_2(f)&=h^{0}(\mb{P}^1,j_{*}\mathcal{L}^{m_1+m_2-2})+\textup{dim}\left(\bigoplus_{x\in \Delta(f_1)}H^{m_2-2}(Y_2)\right)\\
            &=h^{m_1+m_2-2}(Y_1\times Y_2)-h^{m_1-1}(Y_1)h^{m_2-1}_{\textup{van}}(Y_2)+d_1 h^{m_2-2}(Y_2),
        \end{split}
    \end{equation*}
\begin{equation*}
        \begin{split}
            \textup{dim}_{\C}E^{1,m_1+m_2-3}_2(f)&=h^{1}(\mb{P}^1,j_{*}\mathcal{L}^{m_1+m_2-3})\\
            &=d_2h^{m_1-2}(Y_1)-2h^{m_1-2}(Y_1)h^{m_2-1}_{\textup{van}}(Y_2),
        \end{split}
    \end{equation*}
    
and \begin{equation*}
        \begin{split}
            \textup{dim}_{\C}E^{2,m_1+m_2-4}_2(f)&=h^{0}(\mb{P}^1,j_{*}\mathcal{L}^{m_1+m_2-4})\\
            &=h^{m_1+m_2-4}(Y_1\times Y_2)-h^{m_1-3}(Y_1)h^{m_2-1}_{\textup{van}}(Y_2).
        \end{split}
    \end{equation*}
 Again, we get the formula in Lemma \ref{dim} (2) by summing the three dimensions. 

\textbf{Case III$'$:} $H^{m_1-1}_{\textup{van}}(Y_1)=H^{m_2-1}_{\textup{van}}(Y_2)= 0$. 

In this case,  $j^*j_{*}R^kf_{*}\underline{\C}_{\tilde{X}_1\times_{\mb{P}^1}\tilde{X}_2}$ is a locally constant sheaf and the kernel of the adjunction homomorphism $a_k$ is $$\bigoplus_{x\in\Delta(f_1)}H^{m_1-1}(F_x)\otimes H^{k-m_1}(Y_2)\oplus \bigoplus_{x\in\Delta(f_2)}H^{k-m_2}(Y_1)\times H^{m_2-1}(F_x).$$ Hence, \begin{equation*}
        \begin{split}
            \textup{dim}_{\C}E^{0,m_1+m_2-2}_2(f)&=h^{m_1+m_2-2}(Y_1\times Y_2)+d_1 h^{m_2-2}(Y_2)+d_2 h^{m_1-2}(Y_1),
        \end{split}
    \end{equation*}
\begin{equation*}
            \textup{dim}_{\C}E^{1,m_1+m_2-3}_2(f)=0,
    \end{equation*}
    
and \begin{equation*}
            \textup{dim}_{\C}E^{2,m_1+m_2-4}_2(f)=h^{0}(\mb{P}^1,j_{*}\mathcal{L}^{m_1+m_2-4})=h^{m_1+m_2-4}(Y_1\times Y_2).
    \end{equation*}
    
 The formula in Lemma \ref{dim} (2) follows from summing the three dimensions. 
   
\end{proof}

\begin{remark}
    Apply a similar computation, we can prove that the difference of $$ h^{m_1+m_2-1}(\tilde{X}_1\times_{\mb{P}^1}\tilde{X_2})$$ and  the dimension of the degree $m_1+m_2-1$ part of $$H^{\bullet}(\tilde{X}_1)\otimes_{H^{\bullet}(\mb{P}^1)}H^{\bullet}(\tilde{X}_1)$$ is 
    $$\left(h^{m_1}(\tilde{X}_1)-2h^{m_1-2}(X_1)\right)h^{m_2-1}_{\textup{van}}(Y_2)+h^{m_1-1}_{\textup{van}}(Y_1)\left(h^{m_2}(\tilde{X}_2)-2h^{m_2-2}(X_2)\right)$$
 $$+h^{m_1-2}_{\textup{van}}(B_1)\left(h^{m_2-1}(X_2)-h^{m_2-3}(X_2)\right)+\left(h^{m_1-1}(X_1)-h^{m_1-3}(X_1)\right)h^{m_2-2}_{\textup{van}}(B_2)$$
 $$+2h^{m_1-1}_{\textup{van}}(Y_1)h^{m_2-1}_{\textup{van}}(Y_2).$$ Therefore, the difference can be nonzero even if both  $Y_1$ and $Y_2$ have no vanishing cohomology. 
      
\end{remark}

\section{Hodge Numbers of Schoen's Calabi-Yau Manifolds} \label{HNS}

As we mentioned in the introduction, Namikawa and Grassi-Morisson showed that if $E_i\xrightarrow{f_i} \mathbb{P}^1$ are two general rational elliptic surfaces, the Picard group of the fibered product is given by $$\textup{Pic}(E_1\times_{\mathbb{P}^1} E_2)\cong \frac{\textup{Pic}(E_1)\oplus\textup{Pic}(E_2)}{\langle f_1^*\mathcal{O}_{\mathbb{P}^1}(1)\sim f_2^*\mathcal{O}_{\mathbb{P}^1}(1)\rangle}.$$ As a consequence, combined with the fact that the Euler number $e(E_1\times_{\mathbb{P}^1} E_2)$ is zero, it follows that the (upper) Hodge diamond of a Schoen's Calabi-Yau threefold is given by 

$$\begin{array}{ccccccc}
     &&&1&&&  \\
     &&0&&0&& \\
     &0&&19&&0& \\
     1&&19&&19&&1
.
\end{array}$$

An application of Corollary \ref{hformula1}, combined with standard computation of Hirzebruch $\chi_y$ genus, allows us to compute the Hodge diamond of higher-dimensional Schoen's Calabi-Yau manifolds.

 Recall that given any smooth projective variety $X$ of dimension $n$, the Hirzebruch genus $\chi_y(X)$ is defined as the polynomial $$\chi_y(X):=\sum_{p=0}^{n}\chi(X,\Omega^p_X)y^p,$$ where $\chi(X,\Omega^p_X):=\sum_{q=0}^{n}(-1)^q h^q(X,\Omega^p_X)$ is the Euler characteristic of $\Omega^p_X$. By the Hirzebruch-Riemann-Roch theorem, $$\chi(X,\Omega^p_X)=\int_X \textup{ch}(\Omega_X^p)\cdot \textup{td}(X).$$ Therefore, the Hirzebruch genus $\chi_y(X)$ can be expressed in terms of the Chern roots $\gamma_i$ of the tangent bundle $T_X$: $$\chi_y(X)=\int_X \prod_{i}(1+ye^{-\gamma_i})\left(\frac{\gamma_i}{1-e^{-\gamma_i}}\right).$$ Each coefficient of $\chi_y(X)$ is a symmetric polynomial of the Chern roots $\gamma_i$, and hence is a function of Chern classes $c_i(X)$ of $X$. If $X$ is a smooth complete intersection of (not necessarily ample) hypersurfaces of a product of projective spaces, the Chern classes can be computed explicitly by the Euler sequence of the projective space and the normal bundle sequence. Hence, $\chi_y(X)$ can also be determined. In practice, the above process can be easily carried out on a computer. For example, Sage has provided a package to express symmetric polynomials in terms of elementary symmetric polynomials. 
 
 Now, we are ready to compute the Hodge numbers of some examples of higher-dimensional Schoen's Calabi-Yau manifolds.
\begin{example}
    Let $E\rightarrow \mb{P}^1$ be a general rational elliptic surface and $S\rightarrow \mb{P}^1$ be a general Lefschetz fibration of quartic $K3$ surfaces. Consider $X:= E\times_{\mb{P}^1} S$, which is a Schoen's Calabi-Yau fourfold. Since $X$ is a smooth complete intersection of hypersurfaces of $\mb{P}^2\times \mb{P}^1 \times \mb{P}^3$ of degree $(3,1,0)$ and $(0,1,4)$, we can compute that $$\chi_y(X)=2y^4-40y^3+204y^2-40y+2.$$ 
    
    Note also that $S$ is the blowup of a $(4,4)$-complete intersection curve on $\mb{P}^3$. Hence, $h^{2,1}(S)=33$. Then an easy calculation using Corollary \ref{hformula1} shows that the (upper) Hodge diamond of $X$ is $$\begin{array}{ccccccccc}
     &&&&1&&&&  \\
     &&&0&&0&&& \\
     &&0&&11&&0&& \\
     &0&&33&&33&&0& \\
     1&&62&&270&&62&&1 .
\end{array}$$
\end{example}

\begin{example}
Let $S_i\rightarrow \mb{P}^1$ be a general Lefschetz fibration of quartic $K3$ surfaces for $i=1,2$. Consider $X:= S_1\times_{\mb{P}^1} S_2$, which is a Schoen's Calabi-Yau $5$-fold. We can compute that $$\chi_y(X)=-168y^4+1848y^3-1848y^2+168y.$$ 
Then a calculation using Corollary \ref{hformula1} shows that the (upper) Hodge diamond of $X$ is $$\begin{array}{ccccccccccccc}
     &&&&&&1&&&&&&  \\
     &&&&&0&&0&&&&& \\
     &&&&0&&3&&0&&&& \\
     &&&0&&66&&66&&0&&& \\
     &&0&&0&&5&&0&&0&& \\
     &1&&105&&1787&&1787&&105&&1&.
\end{array}$$
\end{example}

\begin{example}
Let $E\rightarrow \mb{P}^1$ be a general rational elliptic surface and $Q\rightarrow \mb{P}^1$ be a general Lefschetz fibration of quintic threefolds. Consider $X:= E\times_{\mb{P}^1} Q$, which is a Schoen's Calabi-Yau $5$-fold. We can compute that $$\chi_y(X)=-100y^4+1100y^3-1100y^2+100y.$$ On the other hand, the discriminant locus $B$ of $Q\rightarrow \mb{P}^4$ is a surface of $(5,5)$-complete intersection. Hence, $h^{2,0}(B)=h^{3,1}(Q)=124$, $h^{1,1}(B)=625$ and $h^{2,2}(Q)=626$.

Then a calculation using Corollary \ref{hformula1} shows that the (upper) Hodge diamond of $X$ is $$\begin{array}{ccccccccccccc}
     &&&&&&1&&&&&&  \\
     &&&&&0&&0&&&&& \\
     &&&&0&&11&&0&&&& \\
     &&&0&&0&&0&&0&&& \\
     &&0&&124&&637&&124&&0&& \\
     &1&&35&&1861&&1861&&35&&1&.
\end{array}$$
\end{example}

\begin{example}
    
Let $S\rightarrow \mb{P}^1$ be a general Lefschetz fibration of quartic $K3$ surfaces and $Q\rightarrow \mb{P}^1$ be a general Lefschetz fibration of quintic threefolds. Consider $X:= S\times_{\mb{P}^1} Q$, which is a Schoen's Calabi-Yau $6$-fold. We can compute that $$\chi_y(X)=2y^6 - 372y^5 + 9174y^4 - 23624y^3 + 9174y^2 - 372y + 2.$$ 
Then a calculation using Corollary \ref{hformula1} shows that the (upper) Hodge diamond of $X$ is $$\begin{array}{ccccccccccccccccc}
     &&&&&&&&1&&&&&&&&  \\
     &&&&&&&0&&0&&&&&&&\\
     &&&&&&0&&3&&0&&&&&&\\
     &&&&&0&&33&&33&&0&&&&& \\
     &&&&0&&124&&629&&124&&0&&&& \\
     &&&0&&0&&33&&33&&0&&0&&&\\
     &&1&&278&&8611&&23442&&8611&&278&&1&&.
\end{array}$$
\end{example}

\setlength\bibitemsep{\baselineskip}  

\printbibliography

\end{document}